# On known and less known relations of Leonhard Euler with Poland[1]

Roman Sznajder[2] (Bowie State University)
e-mail: rsznajder@bowiestate.edu

## I. Introduction

Leonhard Euler (1707-1783) was the central figure in science of the $18^{th}$ century. For almost 60 years, Euler was working practically in all branches of mathematics and mechanics, in addition to undertaking research in astronomy, physics and engineering. All in all, Euler was responsible for about a third of mathematical achievements of his time. Not only was he the greatest mathematician of his era, but also the leading figure behind the reorganization of the research programs of two great academies: the St. Petersburg Academy of Sciences—called in Russia the Imperial Academy of Sciences—and the Berlin Academy of Sciences. He had numerous, but not well known, professional connections with Poland, which are the subject of this work.

Euler's academic activities are divided into three distinct periods:

1. St. Petersburg: 1727-1741
2. Berlin: 1741-1766
3. St. Petersburg: 1766-1783.

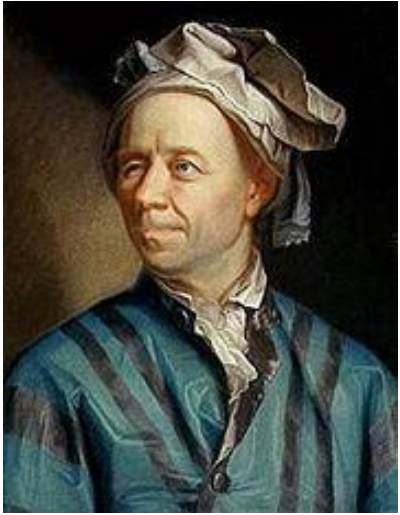 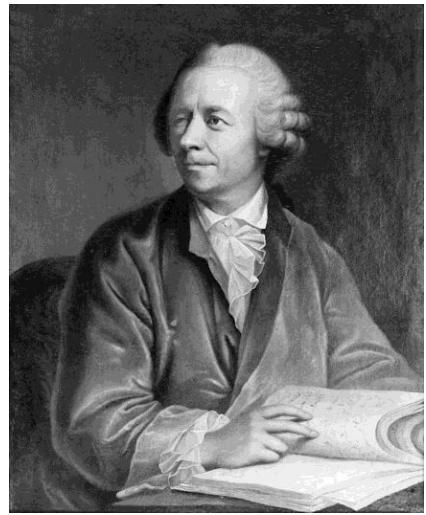

Portraits of Euler by J.E. Handmann[3], 1753 and 1756

---

[1] E-mail address: *rsznajder@bowiestate.edu*. This is an extended version of the presentation given at HPM Americas, March 14-15, 2015, American University, Washington, D.C.
[2] This paper is dedicated to Ms. Anna Parczewska
[3] http://commons.wikimedia.org/wiki/File:Leonhard_Euler.jpg

The first period of Euler's work took place at the St. Petersburg Academy, where he was invited at the initiative of Daniel and Nicolas II Bernoulli, two sons of his teacher, Johann Bernoulli. The Academy was founded by Tsar Peter the Great in 1725. Shortly thereafter, in 1727, Euler arrived in St. Petersburg. In January 1731, he got a professorial position in the department of physics, and in the summer of 1733, he moved to the department of mathematics, where he replaced Daniel Bernoulli who returned to Switzerland. Euler spent 14 years in St. Petersburg. During that period he published 55 volumes, including two-volume *Mechanics* (1736). He left St. Petersburg for Berlin at the invitation of the King of Prussia, Frederick II, whose desire and ambition was to increase the prestige of the Berlin Academy of Sciences.

While leading the mathematics division of the Berlin Academy for 25 years, Euler retained both the honorary membership of the St. Petersburg Academy and his pension. He was publishing in both academies. In addition to his regular obligations, he served the Prussian government as an expert on civil engineering, ballistics, lottery organization, etc. He was also a research editor of the work of both academies, corresponded extensively with external scientific world, and was responsible for hiring new collaborators, purchasing books and research instruments.

## II. Leonhard Euler and Stanisław II August Poniatowski, the Last King of Poland

Euler's stay in Berlin became an opportunity to make acquaintance of several key figures of Polish political and scientific life, most notably, the future king of Poland. At the beginning of 1750, the future—and last— king of Poland, Stanisław II August Poniatowski (1732-1798), visited Berlin for several months for medical treatment. On April 25, 1750, Poniatowski visited the Berlin Academy and most likely met Euler for the first time. Some years later, on November 8, 1764, Euler and his son Johann Albrecht witnessed the reading of the document, entitled *Praise of Stanisław August, the new King of Poland*. Years later, in 1791, Poniatowski also became a member of the Academy [2]. There were several other prominent people from Polish nobility visiting the Academy and having an opportunity to meet Euler.

On June 9, 1766, Prince Adam Kazimierz Czartoryski, in the name of the king, invited Euler with his family to visit Poland. Euler accepted the invitation on the occasion of his trip back to St. Petersburg, where the Imperial Academy had invited him back repeatedly. On his way to Warsaw, he stopped by Poznań and visited a prominent physicist Józef Rogaliński, who built the astronomical observatory on the top of the Jesuit Collegiate in the city. Euler admitted that he "did not expect to see such a quality equipment and mathematical instruments in Poland" [2].

From Poznań, Euler left for Warsaw where he spent ten days as a guest of king Poniatowski. After his stay in Warsaw, Euler left for St. Petersburg through Mitau and Riga, and arrived at his destination on July 28, 1766. Soon after, on August 8, 1766, he sent a letter to the king, warmly thanking him for the hospitality and reception he and his family experienced. Euler's visit led to a decade-long correspondence with the king. The last letter to the king was dated on June 6, 1777, where Euler congratulated him of being elected a member of the St. Petersburg Academy. These letters have been stored in the Archive of Old Documents (the Popiel Collection) in Warsaw. The king followed Euler's work closely, and the cataloguing of the royal collection during 1793-1796 revealed many pieces of Euler's work on various topics.

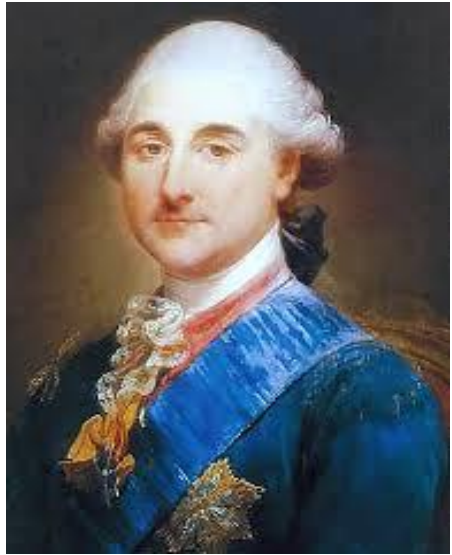

Stanisław II August Poniatowski, king of Poland[4]

Euler's visit in Warsaw went beyond personal contacts between the scientist and the king, and it contributed to Euler's research. For example, using data obtained from the observation of the sun eclipse in 1775 by Jowin F. Bystrzycki, Euler calculated the longitude of Warsaw. More importantly, Euler helped foster collaboration between the Polish and Russian scientific communities. During his short stay in Warsaw, Euler became acquainted with the royal cartographer, Herman Perthees (1739-1814), who sought, by the intercession of the king, information from scientists in St. Petersburg that was relevant to creating maps of Poland and the neighboring states. Ten years later, in his 1776 letter, Euler declared the willingness of the St. Petersburg Academy to cooperate with the royal cartographer on exchange of information and cartography pieces [2].

Euler also had indirect influence on the advancement of Polish education. One of his students, Louis Bertrand (1731-1812) became a teacher of the Swiss mathematician Simon L'Huillier (1750-1840). During 1777-1787, Huillier became a tutor for the aristocratic Czartoryski family. He wrote several textbooks for Polish schools: *Arithmetic for public school*s (1778), *Algebra for public schools* (1782), and *Geometry for public schools* (1780). In the book *Polygonometry* (Gèneve 1789), he continued the work of Euler on polyhedral sets and started a new branch of geometry, *polyedrometry* (science of polyhedrals) [2].

### III. Leonhard Euler and the Danzig Scientific Community

In addition to his ties with the scientific community in Warsaw and Poznań, Euler established close connections with the scientists of the city of Gdańsk (Danzig at the time). Even though Euler never visited Gdańsk, he left a deep mark on its scientific life.

Gdańsk which is more than 1,000 years old, was then a wealthy former member of the Hanseatic League. Even though it had the status of a state-city—also called the Republic of Danzig—it recognized the protection of the king of Poland. The scientific community in Danzig attempted several times to establish its own academic organization, but these efforts failed. Finally, in 1742, the leadership of physicist Daniel Gralath (1708-1767) resulted in the formation of the Danzig Research Society. The latter became the second scientific society in Poland (after the Learned Society, 1720-1727) and the eleventh in the world. Since 1753, the Society also used the name Naturalist Society, but its name in Latin was always *Societas Physicae Experimentalis* [7]. The Society was active until 1945 [8].

---

[4] By Marcello Bacciarelli (1731-1818), ca. 1780; http://niezlomni.com/?p=15243

During his first extended stay in Russia, Euler established ties with several Gdańsk scientists, including the mathematician and astronomer Carl Gottlieb Ehler (1685-1753) who corresponded with Leibniz and became the future mayor of Danzig (1740-1753), and Heinrich Kühn (1690-1769), a professor of mathematics at the Danzig Academic Gymnasium. Two other prominent scientists from Danzig, Johann Phillip Breyne (1680-1764) and Jacob Theodor Klein (1685-1759), started their collaboration with the Imperial Academy [9].

**a. Euler Meets Carl Gottlieb Ehler**

Carl G. Ehler played a dual role of a scientist and diplomat. He was an important member of the delegation of six Gdańsk councilmen who arrived on September 29, 1734, at the court of the Empress of Russia, Anna Ivanovna, to seek forgiveness and reduction of reparations imposed on the city in the aftermath of the War of Polish Succession [9]. On July 9, 1734, after a prolonged siege during this war, Danzig capitulated and was briefly occupied by the Russian army. The city, which supported Stanisław Leszczyński, the losing candidate to the Polish throne, was forced to pay reparations of two million Danzig talars, and the delegation's intent was to obtain forgiveness of the second million. Tough negotiations ensued for several months, but the delegation succeeded and on April 29, 1736, Anna Ivanovna issued the document (*diploma amnestiae*) in which she pardoned Danzig, forgiving the second million of talars and reinstating city privileges [1]. During the duration of negotiations, members of the delegation met with various Russian dignitaries and foreign diplomats. For example, they visited the St. Petersburg Academy of Sciences, Admiralty, shipyards, and metallurgical plants in the vicinity of St. Petersburg. The stay at St. Petersburg gave Ehler the opportunity to meet Euler, either at the end of 1734 or at the beginning of 1735. After the delegation's nearly eight-month-long stay, on May 27, 1735, Empress Anna organized a farewell audience for the delegation. The group left St. Petersburg on June 3, 1735. The details of this visit are included in the report by Carl Ludwig Ehler, the son Carl G. Ehler. The report is deposited in the Gdańsk Scientific Library of the Polish Academy of Sciences (Ms 122) [2].

During his stay in St. Petersburg, Carl G. Ehler participated in numerous meetings of the Imperial Academy, where, on March 7, 1735, Euler presented his *Mechanics*. There was also another scientific connection between Gdańsk and the Imperial Academy as it displayed the letters of one of its most illustrious citizens, Jan Heweliusz (1611-1687)—the most prominent Polish astronomer after Nicolaus Copernicus— on May 13, 1735. These letters were acquired by a French scientist Joseph Nicolas Delisle from the inheritors of Jan Heweliusz in 1726 while Delisle was on the trip from Paris to St. Petersburg [2].

After his return to Gdańsk, Carl G. Ehler corresponded with Euler. The archives of the Russian Academy of Sciences are in possession of fifteen letters from Ehler to Euler and six letters from Euler to Ehler during 1735-1742. On July 15, 1735, Ehler asked Euler for help in recommending his friend and protégé, Heinrich Kühn, to G.W. Krafft in the physics department of the Academy.

In the letter to Krafft of September 24, 1739, Kühn informed of his research interests, including problems in mechanics based on the laws discovered by Christian Wolff, Isaac Newton and Johann Hermann, as well as *Mechanics* by Leonhard Euler. Kühn expressed his wish to became an honorary member of the Imperial Academy. This honor was bestowed upon him on June 27, 1735, and Kühn was granted 100 rubles of annual salary [9].

On May 10, 1740, Euler sent a congratulatory letter to Ehler on the occasion of his election as the mayor of Gdańsk [4], p.387.

**b. Euler, Heinrich Kühn and Other Luminaries**

Heinrich Kühn was born in Königsberg in 1690, where he studied at the Pedagogicum. He then moved to Halle, where he studied law and obtained his JD degree. In 1717, he moved back to Königsberg to continue his studies in natural sciences. He settled in Gdańsk in 1733, where he taught at the Danzig Academic Gymnasium. As a native of Königsberg, Kühn was familiar with the Königsberg Bridges Problem, which asked whether it was possible to design a tour so that one crossed all the bridges over the Pregel river in Königsberg only once.

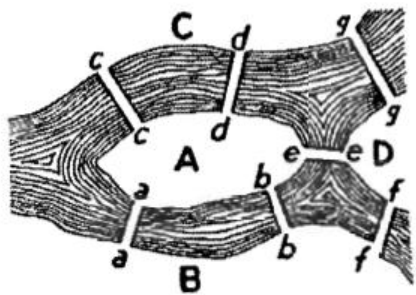

The Königsberg Bridges[5]

It is not completely clear how Euler learned about the Königsberg bridges problem. It is very probable that Carl G. Ehler first discussed this problem with him around 1734/35, during his visit with the Gdańsk delegation in St. Petersburg. As it was already noted by a distinguished Russian science historian, Judith Ch. Kopelevich [9], in his letter of March 9, 1736, Ehler wrote about this problem: "we were discussing it in St. Petersburg." In the same letter, Ehler writes: "You would render to me and our friend Kühn a most valuable service, putting us greatly in your debt, most learned Sir, if you would send us the solution, which you know well, to the problem of the seven Königsberg bridges, together with a proof. It would prove to be an outstanding example of *Calculi Situs*, worthy of your great genius. I have added a sketch of the said bridges ..."

On April 3, 1736, Euler replied: "... Thus you see, most noble Sir, how this type of solution bears little relationship to mathematics, and I do not understand why you expect a mathematician to produce it, rather than anyone else, for the solution is based on reason alone, and its discovery does not depend on any mathematical principle. Because of this, I do not know why even questions which bear so little relationship to mathematics are solved more quickly by mathematicians than by others. In the meantime, most noble Sir, you have assigned this question to the geometry of position, but I am ignorant as to what this new discipline involves, and as to which types of problem Leibniz and Wolff expected to see expressed in this way ..." [4], [5].

---
[5] M. Kraitchik, Mathematical Recreations, Dover, 1953.

Even though, as the quote above indicates, Euler was initially skeptical of the new subdiscipline called the geometry of position (*calculi situs*, mentioned in Ehler's letter), he changed his opinion under the influence of Kühn. On August 26, 1735, Euler presented the solution to the Königsberg bridges problem (in the negative) and its generalizations to the St. Petersburg Academy. It was published in 1736 [3].[6] Following Kühn's suggestion, he incorporated the phrase *geometry of position* into the title of this work. With Euler's paper, graph theory and topology (*geometria situs*) were born. Thus, one can argue that geometria situs indirectly started in Gdańsk [11], [13]. Interestingly enough, Euler never visited Gdańsk or Königsberg.

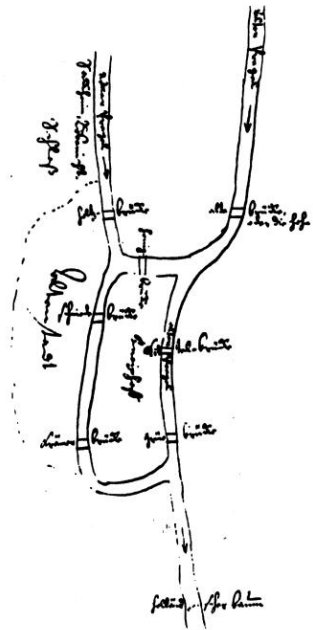

From Ehler's Letter to Euler (1736)

Euler initially corresponded with Kühn through Carl G. Ehler. In the course of their correspondence, Kühn delivered a solution to one of the problems posed by Euler and then, in May of 1735, he sent a manuscript of his paper on some properties of complex numbers, a theme that Euler was very much interested in. In his numerous publications, Kühn presented many original ideas. His crowning achievement was the *Considerations over Constructing Imaginary Quantities and Extracting Imaginary Roots* from 1750/51, published 5 years later in *New Commentaries of the Petersburg Science Academy*, which won him an honor of a corresponding member of the Academy. In addition, another little-known fact is that Kühn was the first person to give a geometric interpretation of complex numbers.

---

[6] For a modern and short proof of the bridges problem, see [14].

In [9] Judith Ch. Kopelevich noted that the correspondence between Euler, Ehler and Kühn represented an interesting contribution to the psychology of scientific research involving team work in solving mathematical problems.

In 1737, Kühn was already directly corresponding with Euler. In his letters, Euler expressed high opinion about Kühn's achievements. There exist 22 letters from Kühn to Euler for the period of 1737-1754, and 2 letters from Euler to Kühn. Euler, through his prolific correspondence, contributed to spreading Kühn's name in the European scientific community. Kühn's paper, *On the origins of water springs and ground water*, which he submitted to the St. Petersburg Academy, drew the attention of numerous scientists. G.W. Krafft and C. Goldbach shared their opinions with Euler. This work, which was a result of Kühn's research in mechanics and hydrology, was earlier awarded a prize by the Bordeaux Scientific Society in 1741 at the competition *Meditationes de origine fontium et atque putealis (On observations of movements of rivers and seas).* The paper was important because it examined the problem of the shape of the earth [2].

As a professor of mathematics at the Danzig Academic Gymnasium, during the years 1735–1770 Kühn became the editor of the calendar, which was distinguished by its high editorial level and interesting content. He became one of the founding fathers of the Danzig Research Society and one of its most active members. In the first of the five treaties printed in 1747 in the *Experiments and Dissertations of the Society*, Kühn described the prototype of an analytical scale and pioneered the theory of scales and weighing [10]. In 1758, on the occasion of the $200^{th}$ anniversary of the Academic Gymnasium, he gave a brilliant and inspiring lecture, *About the Influence of Mathematics and Natural Sciences on the Worldly Happiness of Humankind* [7]. By all means, Heinrich Kühn was the most outstanding Polish mathematician of the $18^{th}$ century. In his letter dated October 28, 1741, Kühn expressed hope that Euler would pay a visit to Gdańsk. Unfortunately, Euler left St. Petersburg for Prussia by ship and stopped only in the city of Szczecin (Stettin).

Even though Euler never made it to Gdańsk, he continued to maintain ties with prominent Gdańsk figures. In 1742, Euler hosted a prominent Gdańsk lawyer, Johann Friedrich Jacobsen. Jacobsen's diary bears the following inscription, written by Euler in May of 1742 (a quote from *Thyestes,* by L. Seneca):

*Nemo tam divos habuit faventes*
*Crastinium ut posset sibi policeri.*

*Berolini, d. 22 May 1742*                                            *memoria causa p.*
                                                                              *Leonh. Euler.*

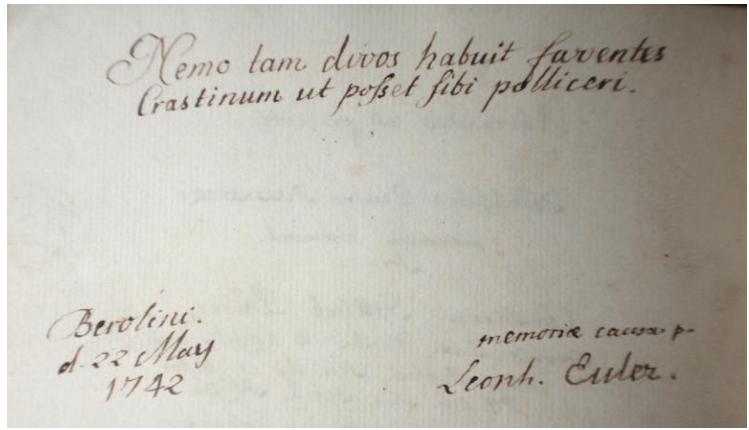

Which means: *No one would have tomorrow such friendly gods as one can expect today.*

*For memory*
*Leonh. Euler.*

The original of this diary [6] is in the Gdańsk Scientific Library of the Polish Academy of Sciences. The diary, in its own way, is an interesting historical artifact. The front page reads:

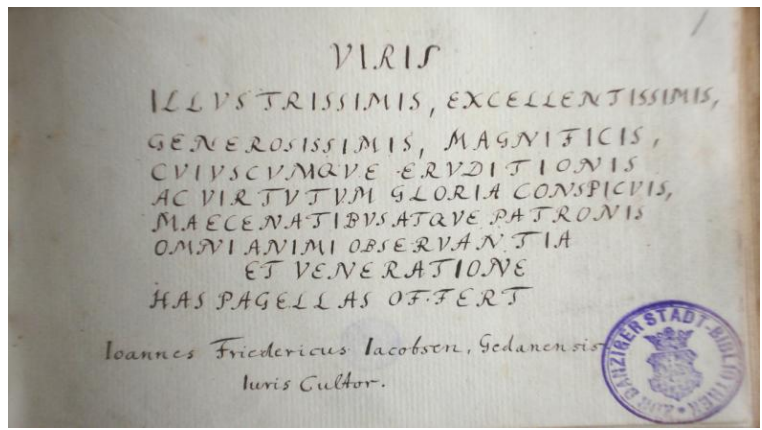

*VIRIS*
*Illustrissimis, excellentissimis,*
*generosissmis, magnificis,*
*cuiuscumque eruditionis*
*acvirtutum gloria conspicuis,*
*maecenatibus atque patronis*
*omni animi observantia*
*et veneratione*
*has pagellas offert*

*Ioannes Fredericus Iacobsen, Gedanensis Iuris Cultor*

This translates as:

*TO GENTLEMEN*
*Most illustrious, most excellent,*
*most generous, most splendid,*
*in visible glory of education and virtues,*
*patrons and curators*
*with wholehearted reverence and admiration*
*this page is being offered by*

*Joannes Fredericus Jacobsen, the Curator of Gdańsk Law*[7]

There were several other Polish scientists who collaborated with the St. Petersburg Academy and had contacts with Euler. Jan M. Hube of Toruń (Thorn), later professor of mathematics and physics, and subsequently the Rector of the famous *Knights' School* in Warsaw, was a student in Gottingen. In his letter of March 12, 1759, he asked Euler to express his opinion on the paper, *Abhandlung von des Kegelschnitten (Treatise on conic sections).* Euler praised the method proposed in the paper in his letter to Hube on April 3, 1759 [2].

**c. Euler and Nathanael Matthaeus von Wolf**

Euler and his family developed particularly close and warm relationship with a Gdańsk physician, astronomer and botanist, Nathanael Matthaeus von Wolf (1724-1784), whom Euler met in 1761. As is apparent from Wolf's letter to the oldest son of Euler, Johann Albrecht, dated March 8, 1768, Wolf and L. Euler became friends in 1766 while the latter was visiting Warsaw. The letter is in the archives of the University Library in Tartu, Estonia.

Nathanael Wolf studied in Jena, Halle, Leipzig and Erfurt. He obtained his degree in medicine in 1768 and became a court doctor for the Lubomirski and Czartoryski Polish aristocratic families. He became a member of the Naturalist Society (1776) and London Royal Society (1777). In 1765-1769, N. Wolf conducted astronomical observations in the Blue Palace in Warsaw while working at the Corps of Cadets as a physician general of the Polish military. In 1768, he was knighted by King Poniatowski. In 1769, Wolf started practicing medicine in Tczew, a town in the vicinity of Gdańsk. In 1772, after the First Partition of Poland, he moved to Gdańsk, which was still under the Polish jurisdiction, so as not to become a Prussian citizen, and opened his doctor's office there. In 1781, using his own funds, Wolf built his astronomical observatory on Bishop's Hill and equipped it with state-of-the-art instruments. Wolf sent the Imperial Academy the results of his observations of the sun eclipse from October 17, 1781. His observatory gained an excellent reputation in the rest of Europe for the accuracy of its observational data. The observatory was destroyed during the Napoleonic war by the Russians during the siege of Gdańsk in 1813.

Among other things, Wolf informed Euler about his botanical work, *Genera plantarum, Vocabulis characteristicis definita* (1776). This is the place where the genus name *Vincetoxicum* was first published. The botanical genus name is derived from the Latin words for *vinci* - win, defeat and *toxicum* for poison. This refers to the purported effects of vegetable juice as an antidote against snake venom[8].

---

[7] I am indebted to Ms. Anna Parczewska for her help with translation of these excerpts.

Wolf was a strong proponent of the inoculation against smallpox, which was plaguing the population every several years. Children were the most vulnerable group. Initially, Gdańsk society was very reserved but gradually acquiesced, and in 1774, the City Council accepted inoculation. N. Wolf was the only doctor who was prepared to perform the procedure.

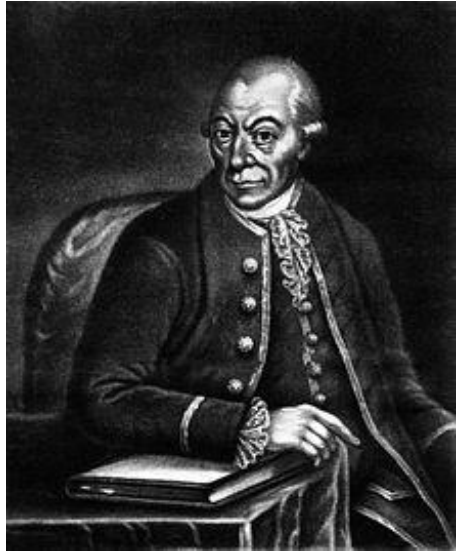

Nathanael Matthaeus von Wolf (1724-1784)[9]

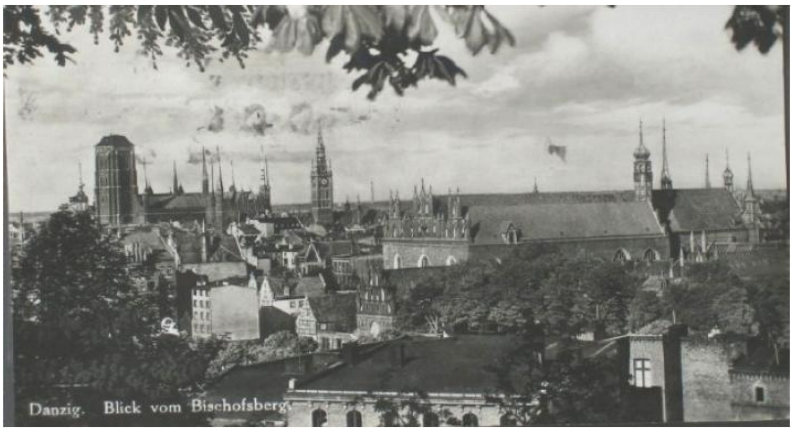

View of the City of Gdańsk from Bishop's Hill[10]

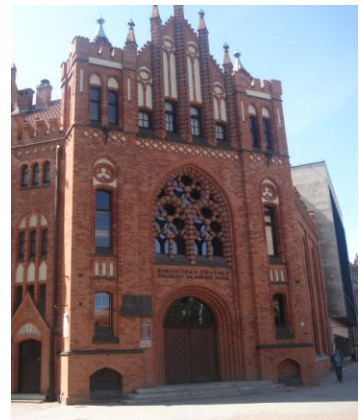

Gdańsk Scientific Library of the Polish Academy of Sciences[11]

---

[8] See http://memim.com/vincetoxicum.html
[9] http://en.wikipedia.org/wiki/Nathanael_Matthaeus_von_Wolf
[10] Pre-WWII postcard of Gdańsk
[11] From the collection of the author

Wolf died in December of 1784 while helping his fellow citizens battle the flu epidemic. His organism, weakened in his youth by tuberculosis, could not defend itself anymore.

**IV. Johann Albrecht Euler and Polish Scientists**

Leonhard Euler shared his passion for science with his oldest son, Johann Albrecht Euler. There are numerous known academic contacts of Johann Albrecht Euler (1734-1800) with Polish scientists. He became the secretary of the Imperial Academy in 1769. In 1794, he became a member of the Danzig Research Society. As the secretary of the Academy, he sent, on December 5, 1792, annals of two periodicals to the Society's library, beginning in the year 1780: "Acta Academiae Scientiarum Imperialis Petropolitanae" and "Nova Acta Academiae..." [2].

Among individual Polish scholars, J.A. Euler corresponded with a Vilnius astronomer Marcin Poczobutt (1728-1810), whom he informed in a letter of December 19, 1777, about the discovery of a new comet by J.A. Lexell.

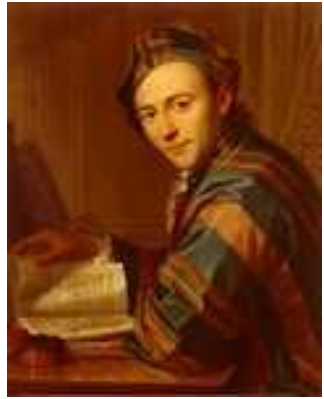

Johann Albrecht Euler (1734-1800)[12]

Another renowned Polish scientist, Jan Śniadecki (1756-1830) also from Vilnius, studied Leonhard Euler's works as a student at the University of Göttingen. Śniadecki also corresponded with Nicolaus Fuss (1755-1826) from Basel, who became a secretary of the Academy in 1800, succeeding J.A. Euler after his death. In 1783, N. Fuss, a former student (and grandson-in-law) of L. Euler and later an editor of his works, wrote a eulogy for L. Euler, which is a masterpiece in its own right. Here is a fragment of Fuss' Eulogy[13] (footnote 8), which highlights Euler-Poland relations:

(8) Throughout his life he preserved the sweet memory of the great goodness that the King showed him and the bonding which inspired his affection towards this heartfelt and spiritual prince and this continued through a correspondence which he had the honor to maintain with him. I cannot withstand the temptation to ornament this eulogy with one of these letters that the king wrote in 1772:

---

[12] By J.E. Handmann (1756), http://www.geni.com/people/Johann-Albrecht-Euler/6000000006713164632
[13] http://www.gap-system.org/~history/Extras/Euler_Fuss_Eulogy.htm

*Monsieur le professor Euler, In response to your letter of 4 August last, I had hoped to be able to confirm the opinion that you are enjoying different circumstances from which your friendship for me has been dictated by the expression of a virtuous and sensible heart. However ... I wish to thank you to good wishes in this regard, and I continue on to the recognition that I owe to your caring in communicating to me the observations of that the accurate astronomers that your Academy made at Bender and near the sources of the Dniestr and the Danube with the locations of some places equally important to geography. I will attempt to put them to good use as with these that are be done in this country with hard work and success, despite the troubles which blocked scientific progress. I request that you continue, as much towards the public utility as my own personal satisfaction and for you to take the opportunity to keep me apprised. I pray to God that he keeps you, dear professor Euler in his holy and worthy safety.*
Written at Warsaw, 7 June 1772

King Stanislaus Augustus                                              (Translation by John S.D. Glaus)

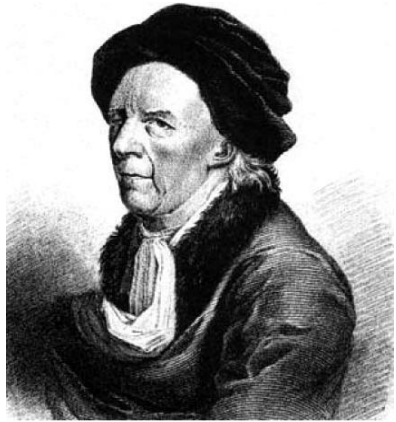

Portrait of Euler by Joseph F.A. Darbes (1778)[14]

**Acknowledgements.** I would like to express my gratitude to Ms. Małgorzata Czerniakowska, a curator of the Gdańsk Scientific Library of the Polish Academy of Sciences, for a long conversation and sharing her work [2]. I am also indebted to the Director of the Library, Dr. Zofia Tylewska-Ostrowska, for giving me permission to photograph pages of J.F. Jacobsen's diary [6] and present in this paper. Finally, I would like to thank an exceptional knot theorist Professor Józef H. Przytycki of The George Washington University, for his suggestions, notes, and a copy of his book [12].

---

[14] https://www.pinterest.com/eulerstalker/leonhard-euler/